\documentclass[10pt,a4paper,twoside]{article}
\usepackage{layout}
\usepackage{exscale, amsmath, amsfonts, amssymb, amsthm, amscd}
\usepackage{verbatim}
\newcommand{\pfskip}{\vspace{3mm}}

\newcommand{\ov}[1]{{\overline #1}}
\newcommand{\wt}[1]{{\widetilde #1}}

\newcommand{\convinfty}[1]
  {\stackrel{#1\to\infty}{-\hspace{-2mm}
  -\hspace{-2mm}-\hspace{-2mm}\longrightarrow}}

\newcommand{\N}{\mathbb{N}}

\newcommand{\R}{\mathbb{R}}
\newcommand{\Z}{\mathbb{Z}}
\newcommand{\T}{\mathbb{T}}

\newcommand{\w}{\omega}
\newcommand{\W}{\Omega}
\newcommand{\MW}{M\times\W}

\newcommand{\B}{{\mathcal B}}

\newcommand{\F}{{\mathcal F}}

\newcommand{\J}{{\mathcal J}}
\renewcommand{\L}{{\mathcal L}}

\newcommand{\Tail}{{\mathcal T}}
\newcommand{\eps}{\varepsilon}

\renewcommand{\t}{\tau}
\newcommand{\tl}{\tau_\lambda}
\newcommand{\ErgSum}[2]{\sum_{i=0}^{n-1}#1\circ #2^i}
\newcommand{\ErgAv}[2]{\frac{1}{n} \sum_{i=0}^{n-1}#1\circ #2^i}

\newcommand{\md}[1]{#1\mbox{\ \text{{\rm mod}}\ }1}

\newcommand{\ArAverage}{\frac{1}{n}\sum_{i=0}^{n-1}}

\newtheorem{thm}{{\bf  Theorem}}[section]
\newtheorem{lemma}[thm]{{\bf  Lemma}}
\newtheorem{cor}[thm]{{\bf  Corollary}}
\newtheorem{prop}[thm]{{\bf  Proposition}}
\newtheorem{definition}[thm]{{\bf  Definition}}

\newtheorem{ex}[thm]{{\bf  Example}}
\newcommand{\pf}{\noindent {\bf Proof. }}

\newcommand{\nn}{\nonumber}

\begin{document}

\pagestyle{empty}

\begin{center}
{\sc
{\large\bf On uniform convergence in ergodic theorems}\\[1mm]
{\large\bf for a class of skew product transformations}
}

\smallskip
{\sc J. Brettschneider}

\end{center}


\begin{abstract} 
\noindent Consider a class of skew product transformations consisting of  
an ergodic or a periodic transformation on a probability space $(M, \B,\mu)$ in the 
base and a semigroup of transformations on another probability space ($\W,\F,P)$ 
in the fibre.  Under suitable mixing conditions for the fibre transformation,
we show that the properties ergodicity, weakly mixing, and strongly mixing 
are passed on from the base transformation to the skew product (with respect to 
the product measure).  We derive ergodic theorems with respect to the skew product  
on the product space.  

The main aim of this paper is to establish uniform convergence with respect
to the base variable for the series of ergodic averages of a function $F$ 
on $M\times\W$ along the orbits of such a skew product. 
Assuming a certain growth condition for the coupling function,
a strong mixing condition on the fibre transformation, and continuity and
integrability conditions for $F,$ we prove uniform convergence in the base and 
$\L^p(P)$-convergence in the fibre.  Under an equicontinuity assumption on $F$ 
we further show $P$-almost sure convergence in the fibre.
Our work has an application in information theory:   
It implies convergence of the averages of functions on random fields
restricted to parts of stair climbing patterns defined by a direction.
\end{abstract}

\renewcommand{\thefootnote}{}
\vspace{-5mm}

\section{Introduction}\label{int}

The approximation of a line by a planar lattice yields a stair climbing pattern.
Consider the averages of a function of a random field along a finite window
moving up the stair climbing pattern.  
Under which conditions does this sequence converge, and what are explicit
formulae for the limit?
More formally, let $L_{\lambda,t}(z):=(z,[\lambda z+t])\ (z\in \Z)$ be a lattice 
approximation of  the line with slope $\lambda$ and $y$-intercept $t,$
and let $m\in\N$ be a fixed window size. 
Let $P$ be a $\Z^{2}$-indexed random field with values in a set
$\Upsilon,$ i.e., a stationary probability measure on $\W:=\Upsilon^{{\Z^{2}}},$
and let $\F$ be the canonical $\sigma$-algebra.  
Consider the averages
\begin{equation}\label{Intro:ex.sq}
   \frac{1}{n}\sum_{i=0}^{n-1}f \big( \w(L_{\lambda,t}(i,...,I+m-1)) \big)
\qquad(n\in\N)
\end{equation}
of a function $f\in \L^{1}(\W,\F,P).$
What can we say about $P$-almost sure or $\L^{1}(P)$-convergence
of this sequence?  
Averages of this type are similar to the ones that occur in the context of directional
Shannon-MacMillan theorems for lattice random fields (cf.~\cite{Bre01SMcM}).  

The situation described above can be represented as a special case of a 
more general set up involving skew product transformations.
In the independent component, or {\it base,} we have a measure-preserving
transformation $\t$ on a probability space $(M,\B,\mu).$   In the dependent 
component, or {\it fibre,} we have a mixing semigroup of measure-preserving 
transformations $(\theta_{k})_{k\in K}$ on a probability space $(\W,\F,P),$
which is linked to the base by a $K$-valued $\B$-measurable function 
$\kappa$ on $M.$   The class of skew products considered in this paper is
given by
\begin{eqnarray}\label{Intr:def.skew}
  S(t,\w)=(\t (t),\theta_{\kappa(t)}\w)\qquad(t\in M,\w\in\W).
\end{eqnarray}
This paper deals with the following two questions: \\
\noindent (A) Is $S$ ergodic or mixing with respect to the product measure?\\
\noindent (B) Do the ergodic averages along 
the skew product converge uniformly with respect to the base?\\
Under suitable assumptions, we will answer both questions positively.

Ergodicity, and other mixing properties of various classes of skew products 
have been studied by a number of authors, but the above situation does 
not fit into any of the settings covered by existing literature. 
Kakutani \cite{Kak51} introduced a skew product with a Bernoulli-shift 
in the base and an ergodic transformation in the fibre was introduced by . 
He showed that the skew product is ergodic if and
only if the transformation in the fibre is ergodic. Other mixing
properties were investigated, e.g., by Meilijson \cite{Mei74},
den Hollander and Keane \cite{dHK86}, 
and Georgii \cite{Geo97}. Adler and Shields (cf.~\cite{AdS72} and
\cite{AdS74}) considered a translation on the cirle for the fibre.
Anzai \cite{Anz51} introduced skew products of two translations
on the circle, and derived a criterion for ergodicity. Furstenberg
\cite{Fur61} studied unique ergodicity. Zhang \cite{Zha96} investigated
this for a translation on a torus in the fibre.
A torus translations in the base can also be combined with the
translation on $\R$ by the value of a real function of the argument
in the base. Skew products of this type are called real extensions of
torus translations, and they were explored in Oren \cite{Ore83},
Hellekalek and Larcher \cite{HeL86}, \cite{HeL89}, and Pask \cite{Pas90}.

\pagestyle{plain}
  
Our answer to question (A) is summarized in Theorem \ref{ErgS}.  
We prove that, under suitable mixing conditions 
for the transformations in the fibre, the properties ergodic, weakly mixing, 
and strongly mixing, are passed on from the transformation in the base 
to the skew product.  As an explicit example we study the case when $P$ is 
a random field and $(\theta_{k})_{k\in K}$ is a group of shift transformations.
In this case, the conditions on the fibre transformation can be insured 
by assuming tail-triviality for $P$ and a growth condition for ergodic sums of
$\kappa$ along $\t$ (cf.~Corollary \ref{ShiftsSProp}).

Our answer to question (B) is given in Theorem \ref{UnifL}.
The proof combines two different approaches. The first approach explores
uniform convergence theorems in the spirit of Weyl's classical result
for the rotation on the circle.  As a little addition to the theorems of Weyl and 
Oxtoby we show that the ergodic averages of a continuous and uniquely ergodic 
transformation on a real interval convergence uniformly for the class of 
Riemann-intergrable functions  (cf.~Corollary \ref{UnifR}). 
The second approach are techniques developed for ergodic theorems 
along subsequences.  We extend Blum and Hanson's theorem to the
$d$-parameter case replacing the strict monotonicity condition on the sequence
by a growth condition on the coupling function $\kappa,$ uniformly in $t.$
Combining the two approaches we obtain uniform convergence in the base and 
$\L^p$-convergence in the fibre for functions that are continuous with respect
to the base and fulfill an integrability condition in the fibre (cf.~Theorem \ref{UnifL}).
We further derive a variation of this theorem for the class of Riemann-integrable
functions on a real interval (cf.~Corollary \ref{UniLRMak}). 
We further derive a result (cf.~Theorem \ref{UnifP}) about uniform convergence 
in the base and $P$-almost sure convergence in the fibre, provided the 
iterates of the function fulfill an equicontinuity condition. 

We conclude the paper by returning to our initial questions about the asymptotics of 
(\ref{Intro:ex.sq}).  
Let $\T:=[0,1)$ be the circle equipped with the Borel $\sigma$-algebra 
$\B$ and the Lebesgue measure $\mu.$ 
For $\lambda\in\R$ define a rotation on $\T$ by
$\t_\lambda(t):=  t+\lambda\  {\rm mod}\ 1.$ 
For $x\in\R$ let $[x]$ be the integer part of $x.$
Define a skew product on the product space $\T\times\W$ by
\begin{equation}\label{Intro:ex.sp}
S(t,\w):=\big(\t_{\lambda}(t),\vartheta_{(1,[\lambda +t])}\w\big)\qquad (t\in\T,\w\in\W).
\end{equation}
We will see that the iterates of $S$ follow the stair climbing pattern $L_{\lambda,t}.$ 
This allows to rewrite the sequence (\ref{Intro:ex.sq}) as an average of 
$f$ along the orbit of $S.$  The convergence of this sequence, 
for {\it all} starting levels $t$ of the stair climbing pattern, 
is  a consequence of the uniform ergodic theorems derived in in Section \ref{ret}
(cf.~Corollary \ref{Unif:ex.staircase}). 

{\bf Ouline of the paper: }
In the first section we define mixing properties of semigroups of transformations
along sequences and we introduce a the class of skew products considered in this paper.
In Theorem \ref{ErgS}, we give the result on ergodicity and mixing properties of these 
skew products.  Finally, we have a closer  look on the case when the fibre transformation 
is a shift operator for a random field.  
In Section \ref{ret} we discuss ergodic theorems for skew products
with ergodic base transformation (cf.~Corollary \ref{Reterg}) and with 
periodic base transformations (cf.~Corollary \ref{Retper}).  We illustrate the results
with two examples related to the sequence \eqref{Intro:ex.sq}.
In the last section we focus on the main aim of this paper, the uniform
convergence with respect to the base.  Depending, among other things, 
on the regularity of the function with respect to the base variable, we obtain
different kinds of convergence in the fibre.
Theorem \ref{UnifL} states $\L^p(P)$-convergence provided the function is
continuous with respect to the base variable.  Corollary \ref{UniLRMak} is a version 
of this for Riemann-integrable functions.
Theorem \ref{UnifP} states $P$-almost sure convergence, provided the functions
fulfill a certain equicontinuity condition.
Corollary \ref{Unif:ex.staircase} brings us back to the original motivation for this paper.
It states the convergence of the sequences \eqref{Intro:ex.sq}. 

\section{Mixing properties of a class of skew products}\label{spc}

Let $(\W,\F,P)$ be a probability space,
and let $(\theta_k)_{k\in\N_{0}^{d}}$ be a {\it d-parameter semigroup} of
measure-preserving transformations on $(\W,\F,P),$ i.e.,
each of the transformations preserves the measure $P,$ and
$\theta_0=\text{{\rm Id}},$ and 
$\theta_k\circ\theta_l=\theta_{k+l}\ \text{for all}\ k,l\in\N_0^d.$
$(\theta_k)_{k\in\Z_{0}^{d}}$ is a {\it d-parameter group} if
$(\theta_k)_{k\in\N_{0}^{d}}$  is a semigroup and
$\theta_{-k}=\theta_{k}^{-1}$ for all $k\in \N^d.$
The following example will be used frequently in our settings.
Let $\sigma_1$ and $\sigma_2$ be two commuting
measure-preserving transformations on $(\W,\F,P).$ Then
\begin{eqnarray}\label{SegrEx}
  \theta_k :=\sigma_1^{k^{(1)}}\circ\sigma_2^{k^{(2)}}
  \quad\hbox{for}\quad k=(k^{(1)},k^{(2)})\in\N_0^2
\end{eqnarray}
defines a two-parameter semigroup $(\theta_k)_{k\in \N_0^2}$ of
measure-preserving transformations on $(\W,\F,P).$  
If $\sigma_1$ and $\sigma_2$ are invertible
it extends to a two-parameter group $(\theta_k)_{k\in \Z^2}.$
The constructions extends to $d$-parameters in a obvious way.

Let $K=\N_0^d$ or $K=\Z^d.$  Let $\t$ be a measure-preserving 
transformation of a probability space $(M ,\B ,\mu ),$
and assume that $\kappa$ is a $\B$-measurable function 
on $M$ with values in $K.$ Then
\begin{eqnarray*}
   S(t,\w)=(\t (t),\theta_{\kappa(t)}\w)\qquad (t\in M,\w\in\W)
\end{eqnarray*}
defines a skew product on the product space $\ov\W:=M\times\W.$
In particular, choosing $\kappa\equiv k_0$ for a constant $k_0\in K$
yields the uncoupled product of $\t$ and $\theta_{k_0}.$
Obviously, $S$ is measurable with respect to the product $\sigma$-algebra 
$\sigma$-algebra $\ov\F:=\B\otimes\F,$ and it preserves the product measure 
$\ov P:=\mu\otimes P.$
It is easy to see that for all $n,m\in\N_0,$ for all $t\in M,$ and for all $\w\in\W,$
\nopagebreak
\begin{eqnarray}\label{ItS}
 S^n (t,\w) = \big(\t^n (t),\theta_{\kappa_n(t)}\w\big)
\quad \mbox{and}\quad
 \kappa_{n+m}=\kappa_n+\kappa_m\circ\t^n, 
\text{  where }
\kappa_n = \ErgSum{\kappa}{\t}.
\end{eqnarray}

Let us know study under which conditions the skew product is ergodic.
Furthermore, as suggested by J. Aaronson, we broaden the question to other 
mixing properties.  Answers to these questions will be given in the next lemma.  
Note that, by a simple projection argument, the ergodicity of $\t$ is necessary
for the ergodicity of $S.$  As we know in the case of an uncoupled product,
the ergodicity of both transformations does not garanty the ergodicity of
the product.  However, it can be shown that the product
is ergodic whenever one of the transformations is ergodic and the
other one is weakly mixing (cf.~\cite{Kre85}).
A major ingredient for $S$ for our lemma are assumptions that 
bring into play the function $\kappa.$  We make use of two conditions:
\begin{itemize}
\item[{\bf (C1)}] $(\theta_k)_{k\in K}$ is weakly mixing along the sequence
            $(\kappa_n(t))_{n\in\N}$
            for $\mu$-almost all $t\in M.$
\item[{\bf (C2)}] $(\theta_k)_{k\in K}$ is strongly mixing and
            $(\kappa_n(t))_{n\in\N}$ goes to infinity for
            $\mu$-almost all $t\in M.$
\end{itemize}

\noindent Note that (C2) implies (C1). 
Condition (C2) can be easily verified for lattice approximations of a line, 
as discussed in the context of Corollary \ref{Unif:ex.staircase}.  
We further need the notion of
{\it weakly mixing along a sequence.}  This has been introduced for 
transformations by N.\ Friedman (cf.~\cite{Fri83}), and 
we extend it to $d$-parameter (semi-)groups.

\begin{definition}\label{SeqMixGr}  Let $(k_n)_{n\in\N}$ be a $K$-valued 
sequence.  A (semi-)group $(\theta_k)_{k\in K}$ of
measure-preserving transformations on $(\W,\F,P),$
is called {\rm weakly mixing along $(k_n)_{n\in\N}$} with respect to $P,$ if
\begin{eqnarray}\label{knwmixGr}
  \ArAverage\big\vert P\big(A\cap\theta^{-1}_{k_i}B\big)
  - P\big(A\big)\,P\big(B\big)\big\vert\convinfty{n}0\qquad
  \text{for all }A,B\in\F.
\end{eqnarray}
\end{definition}

\begin{lemma}\label{ErgS} 
(i) Assume condition (C1).  If $\t$ is ergodic w.r.t.~$\mu$
then $S$ is ergodic w.r.t.~$\ov{P}.$

\noindent (ii) Assume condition (C1).  If $\t$ is weakly mixing 
w.r.t.~$\mu$ then $S$ is  weakly mixing w.r.t.\ $\ov{P}.$

\noindent (iii) Assume condition (C2).  If $\t$ is strongly mixing 
w.r.t.~$\mu$ then $S$ is strongly mixing w.r.t.~$\ov{P}.$
\end{lemma}
\pf
We give the proof of the first statement here; the remaining proofs are
conducted in a similar fashion.  Assume condition (C1) and
let $\t$ be ergodic with respect to $\mu.$
To prove the ergodicity of $S$ we will show that for all bounded
$\ov\F$-measurable functions $F$ and $G$
\begin{eqnarray*}
\frac{1}{n}\sum_{i=1}^n \int_{\MW}F\circ S^i(t,\w)\cdot
G(t,\w)\,d\ov P \convinfty{n}\wt F\cdot\wt G,
\end{eqnarray*}
where
$$
\wt F =\int_{\MW}F(t,\w)\,d\ov P
\quad\quad\text{and}\quad\quad
\wt G =\int_{\MW}G(t,\w)\,d\ov P.
$$
It is sufficient to show this for functions which are products of
functions on the factors, i.e.,
$F(t,\w)=f(t)\Phi(\w)\text{ and }G(t,\w)=g(t)\Psi(\w),$
where $f$ and $g$ are bounded $\B$-measurable functions on $M$ 
and $\Phi$ and $\Psi$ are bounded $\F$-measurable functions on $\W.$ 
(The general case follows by approximation.)  For these functions we have
$
\wt F=\wt f\cdot\wt\Phi \text{ and }
\wt G=\wt g\cdot\wt\Psi,
$
with
$$
\wt f=\int_M f\,d\mu,\ \ \wt\Phi=\int_\W \Phi\,dP,\ \
\wt g=\int_M g\,d\mu,\ \ \text{and}\ \ \wt\Psi=\int_\W\Psi\,dP,
$$
and we obtain
\begin{align}\nonumber
  \bigg\vert \,\frac{1}{n}&\sum_{i=1}^n \int_{\MW} F\circ S^i(t,\w)\cdot
  G(t,\w)\,d\ov P\, -\, \wt F\cdot\wt G\,\bigg\vert\\
  \nonumber
  \le & 
  \bigg\vert\, \int_M \frac{1}{n}\sum_{i=1}^n  \int_{\MW}
  f(\t^i(t))g(t)\Phi(\theta_{\kappa_i(t)}\w)\Psi(\w)
  -f(\t^i(t))g(t)\wt\Phi\,\wt\Psi
  \, d\ov P\,\bigg\vert\,d\mu\\
  \nonumber
  & +
  \bigg\vert\,\frac{1}{n}\sum_{i=1}^n\bigg(\int_{M}
  f(\t^i(t)) g(t)\,d\mu 
  \,-\,\wt f\cdot\wt g\bigg)\cdot\wt\Phi\cdot\wt\Psi\,\bigg\vert\\
  \nonumber
  \le &
  \|f\|_\infty\|g\|_\infty\, \frac{1}{n}\sum_{i=1}^n\bigg\vert\,
  \int_\W\Phi(\theta_{\kappa_i(t)}\w)\Psi(\w)\,dP
  -\wt\Phi\cdot\wt\Psi\,\bigg\vert
  + 
  \bigg\vert\,\frac{1}{n}\sum_{i=1}^n\int_M f\circ\t^i\hspace{0.05mm}\cdot\hspace{0.05mm}
  g\,d\mu \,-\,\wt f\cdot\wt g\,\bigg\vert\,\vert\,\wt\Phi
  \,\vert\,\vert\,\wt\Psi\,\vert.
\end{align}
This goes to $0,$ because
the averages in the first expression convergence to $0$ 
for $\mu$-almost all $t$ by condition (C1) and
the second expression converges to $0$ because of the ergodicity of $\t.$
\qed
\pfskip

This section concludes with a closer look at the case of shift
transformations on a discrete random field with values in a set $\Upsilon.$
Let $\W:=\Upsilon^{\Z^d}.$ For any $J\subseteq\Z^d$ let 
$\F_J$ denote the $\sigma$-algebra generated
by all projections $\w\mapsto\w(j)$ with $j\in J,$
and let $\F:=\F_{\Z^d}.$  Denote the coordinates of an elements in $\Z^d$ by 
upper indices, and let $\|\,\cdot\,\l$ be its maximum norm.
Consider the {\it shift transformations}
$(\vartheta_v)_{v\in\Z^d}$ on $\W,$ i.e.,
$\vartheta_v(\w)(j):=\w(j+v)$ $(j\in\Z^d).$
Let $P$ be a {\it random field}, i.e., a measure on $(\W,\F)$
which is invariant with respect to $\vartheta_v$ for all
$v\in\Z^d.$ The {\it tail field} is the $\sigma$-algebra
$
\Tail:=\bigcap_{V\subset\Z^d\ {\rm finite}}
\F_{\Z^d\setminus V},
$
and $P$ is called {\it tail-trivial} if it fulfills a $0$-$1$ law on $\Tail.$  
Let $v_1, v_2\in\Z^d.$  As in \eqref{SegrEx}, 
$\theta_k:=\vartheta_{v_1}^{k^{(1)}}\circ \vartheta_{v_2}^{k^{(2)}}$ $(k\in\Z^d)$
defines a 2-parameter group of measure-preserving transformations.  We have
\begin{equation}\label{spc:twoShifts}
\theta_{\kappa_n(t)}=\vartheta_{v_1}^{\kappa_n^{(1)}(t)}\,\circ\,\vartheta_{v_2}^{\kappa_n^{(2)}(t)}
=\vartheta_{
       {\kappa_n^{(1)}(t)\,v_1 } \,+\, 
       {\kappa_n^{(2)}(t)\,v_2 } 
       }.
\end{equation}
In this situation we have the following

\begin{cor}\label{ShiftsSProp}
Let $v_1$ and $v_2$ be linear independent vectors in $\Z^d.$
Assume that $P$ is tail-trivial and that the sequence
$\big(\|\kappa_n(t)\|\big)_{n\in\N}$ goes to infinity for
$\mu$-almost all $t\in M.$
Then, when $\t$ is ergodic, weakly
mixing or strongly mixing with respect to $\mu,$ $S$ is ergodic,
weakly mixing or strongly mixing with repect to $\ov P,$
respectively.
\end{cor}
\pf We are going to show condition (C2). Define the boxes
$V_n=\big\{v\in\Z^d\,\big\vert\,\|\,v\,\|\le n\big\}$ $(n\in\N),$ and
let $B\in\F_J,$ for some finite subset $J$ of $\Z^d.$ Then there
is an $m\in\N$ such that $J\subseteq V_m.$ Setting
$m(n):=\kappa_n^{(1)}(t)\,v_1+\kappa_n^{(2)}(t)\,v_2$ we observe
that the translated sets $J-m(n)$ are contained in
$V_{{\widetilde m}(n)}^c,$ where ${\widetilde m}(n)=(m(n)-2m)\vee
0.$ For any $A\in\F,$ we obtain
\begin{align}\label{ShiftsSProp1}
 \nonumber
 \big\vert  P\big(A\cap\theta_{\kappa_n(t)}^{-1}B\big)
          -P\big(A\big)\,P\big(B\big)\big\vert
& =  \big\vert P\big(A\cap
     \vartheta_{\kappa_n^{(1)}(t)\,v_1+\kappa_n^{(2)}(t)\,v_2}^{-1}B\big)
      -P\big(A\big)\,P\big(B\big)\,\big\vert\\
  \nonumber
 &\le\sup_{C\in\F_{\Z^d\setminus V_{{\widetilde m}(n)}}}
     \big\vert\,P(A\cap C)-P(A)P(C)\,\big\vert.
\end{align}
By the assumptions on $v_1, v_2$ and $\kappa,$
$\|{\widetilde m}(n)\|$ goes to infinity.
By Proposition 7.9 in \cite{Geo88}, tail-triviality is equivalent to {\it short-range correlations,} i.e.,
$
\sup_{C\in\F_{\Z^d\setminus V_n}}
\big\vert\,P(A\cap C)-P(A)P(C)\,\big\vert
$
converges to $0$ as $n$ goes to infinity.
\qed

\section{Ergodic theorems with skew products}\label{ret}

Applying Birkhoff's ergodic theorems to the skew product
$S$ yields, for any $F\in\L^1(\ov \W,\ov \F,\ov P),$
\begin{eqnarray}\label{BirkS}
  \frac{1}{n}\sum_{i=0}^{n-1}F\big(\t^i(t),\theta_{\kappa_i(t)}\w\big)
  \convinfty{n}\ov E[F\vert\J]
  \qquad \ov P\text{-almost surely and in }\L^1(\ov P),
\end{eqnarray}
where $\J$ is the $\sigma$-algebra of all $S$-invariant sets in $\ov\F.$
We study  this limit more closely for two different cases:
when the transformation $\t$ is ergodic and when it is periodic.
In the ergodic case, combining $\eqref{BirkS}$ and Lemma \ref{ErgS} 
immediately yields the following ergodic theorem for the skew product.

\begin{cor}\label{Reterg}  
Assume that $\t$ is ergodic with respect to $\mu$ and that the condition
(C1) is fulfilled.
Then for any function $F\in\L^1(\ov\W,\ov\F,\ov P),$
\begin{eqnarray*}
 \frac{1}{n}\sum_{i=0}^{n-1}F\big(\t^i(t),\theta_{\kappa_i(t)}\w\big)
  \convinfty{n}\ov E[F]
  \qquad \ov P\text{-almost surely and in } \L^1(\ov P).
\end{eqnarray*}
\end{cor}

Now consider the case that $\t$ periodic.
We calculate the iterates of the skew product 
and derive an ergodic theorem with an explicit expression for the limit.

\begin{lemma}\label{ItSper}
Assume that $\t$ is periodic with $q\in\N.$ Then for all $j\in\Z$
and all $\nu\in\{0,1,...,q-1\},$
\begin{itemize}
\item[(i)]$\kappa_{jq+\nu}=j\kappa_{q}+\kappa_{\nu},$
\item[(ii)]$\theta_{\kappa_{jq+\nu}(t)},
           =\Big(\theta_{\kappa_{q}(t)}\Big)^j\circ\theta_{\kappa_{\nu}(t)}
           \qquad\text{for all }t\in M,$
\item[(iii)]
    $S^{jq+\nu}(t,\w)= \big (\t^{\nu}(t),
    \theta_{\kappa_{\nu}(t)}\circ(\theta_{\kappa_q(t)})^j\w\big )
    \qquad\text{for all }t\in M \text{ and all }\w\in\W.$
\end{itemize}
\end{lemma}

\pf
(i) follows from the definition of $\kappa$ using the periodicity of $\t.$  
(ii) is an immediate consequence of (i) and the semigroup property of $\theta.$
(iii) follows from (\ref{ItS}), the periodicity of $\t$ and because of (ii). \qed

\begin{cor}\label{Retper}  
Assume $\t$ is periodic with period $q\in\N,$ and $F\in\L^1(\ov\W,\ov\F,\ov P).$
Denote by $\J_t$ the $\sigma$-algebra of
$\theta_{\kappa_q(t)}$-invariant sets in $\F.$  Then
\begin{displaymath}\label{RetperA}
  \frac{1}{n}\sum_{i=0}^{n-1}F\big(\t^i(t),\theta_{\kappa_i(t)}\w\big)
  \convinfty{n}
  \frac{1}{q}\sum_{\nu=0}^{q-1}E\big[F\big(\t^\nu(t),
  \theta_{{\kappa_\nu}(t)}(\cdot)\big)\big\vert\, \J_t\big],
\end{displaymath}
for $\mu$-almost all
$t\in M,$ and for $P$-almost all $\w\in\W$ and in $\L^1(P).$
If $\theta_{\kappa_q(t)}$ is ergodic with respect to $P$
for $\mu$-almost all $t\in M,$ then the limit
simplifies to the constant
$  1/q\sum_{\nu=0}^{q-1} E\big[F\big(\t^\nu(t),\,\cdot\,\big)\big].$
\end{cor}

\pf 
Any $n\in\N$ can be represented as
$n=mq+\nu,$ with $m\in\N$ and $\nu\in\{0,1,...,q-1\},$
and we may break down the ergodic averages to
\begin{displaymath}
  A_nF:=\ErgAv{F}{S}=\frac{mq}{mq+\nu}\bigg(
  \frac{1}{mq}\sum_{i=0}^{mq-1}F\circ S^i
  +\frac{1}{mq}\sum_{i=mq}^{mq+\nu-1}F\circ S^i\bigg).
\end{displaymath}
Since the first factor converges to $1,$ and the second addend within the
brackets converges to $0,$
our question reduces to the study of ergodic limits along the subsequence
$(mq)_{m\in\N}.$  They take the form
  $A_{mq}F=1/q\sum_{\nu=0}^{q-1}A_m^{(\nu)}F,$ where
\begin{eqnarray*}
  A_m^{(\nu)}F(t,\w)
    := \frac{1}{m}\sum_{j=0}^{m-1}F\circ S^{jq+\nu}(t,\w)
    =\frac{1}{m}\sum_{j=0}^{m-1}
    F\big(\t^{\nu}(t),  \theta^{\kappa_{\nu}(t)}\circ(\theta^{\kappa_q(t)})^{\,j}\w\big).
\end{eqnarray*}
The last equality can be seen by applying Lemma \ref{ItSper}.
For $\mu$-almost all $t\in M,$  
the function $f^{(\nu)}_t(\w):= F\big(\t^\nu(t),\theta_{\kappa_\nu(t)}\w\big)$ $(\w\in\W)$
is integrable, and applying Birkhoff's ergodic theorem yields 
\begin{displaymath}
  \lim_{m\to\infty} A_m^{(\nu)}F(t,\cdot)
  =E[f^{(\nu)}_t\vert\J_t]
  =E\big[F\big(\t^\nu(t),\theta_{\kappa_{\nu}(t)}(\cdot)\big)\vert\J_t]
  \qquad\mbox{$P$-almost surely and in $\L^1(P)$}
\end{displaymath}
This implies the first statement of the Corollary.  In the ergodic case $\J_t$ is trivial,
and, using the invariance of $P$ under $\theta,$ the last expression reduces to 
$E\big[F\big(\t^\nu(t),\,\cdot\,\big)].$ 
\qed\pfskip

We end this section by illustrating the results by two special cases relevant to the
skew product tracing a stair climbing pattern introduced in (\ref{Intro:ex.sp}).  

\begin{ex}\label{TorTransl} Circle rotations as base transformations. 
{\rm
Let $\T:=[0,1)$ be the circle equipped with the Borel $\sigma$-algebra 
$\B$ and the Lebesgue measure $\mu.$ 
For $\lambda\in\R$ define a rotation on $\T$ by
$\t_\lambda(t):=  t+\lambda\  {\rm mod}\ 1.$ 
The rotation preserves the measure $\mu$ and it is continuous.
For rational $\lambda$ the rotation is periodic, for 
irrational $\lambda$ it is uniquely ergodic.
The circle can be equipped with a metric $d(s,t):=|s-t|\;(s,t\in\T).$
The metric is rotation invariant.

\noindent
(i) Let $\lambda$ be irrational.  
Choose $(\theta_k)_{k\in K}$ and $\kappa$ fulfilling condition (C1).  
By Lemma \ref{ErgS}, $S$ is ergodic.
By Corollary \ref{Reterg}, for any integrable function $F$ on
$(\T\times\Omega,\B\otimes\F,\mu\otimes P),$
\begin{eqnarray*}
  \frac{1}{n}\sum_{i=0}^{n-1}F(\md{t+i\lambda},\theta_{\kappa_i(t)}
  \omega)\convinfty{n}\int_0^1 E[F(t,\cdot)]\,dt,
\end{eqnarray*}
for $\mu\otimes P$-almost all $(t,\w)\in\T\times\Omega$
and in $\L^1(\T\times\Omega,\B\otimes\F,\mu\otimes P).$

\noindent
(ii) Let $\lambda$ be rational.  
There is a unique representation
$\lambda=p/q,$ where $p\in\Z, q\in\N,$ $p$ and $q$ have
no common divisor. $\t_\lambda$ is periodic with period $q.$ 
Furthermore, $\t_\lambda$ respects the partition
$[0,\frac{1}{q}), [\frac{1}{q},\frac{2}{q}),...,[\frac{q-1}{q},1)$ of $\T,$ 
i.e., for every $\nu\in\{1,...,q-1\}$ there is a
$\widetilde\nu\in\{1,...,q-1\}$ such that
$\tl([\frac{\nu-1}{q},\frac{\nu}{q}))=[\frac{\widetilde\nu-1}{q},\frac{\widetilde\nu}{q}).$
The limit in Corollary \ref{Retper} is of the form
 $  1/q\sum_{\nu=0}^{q-1}E\big[F\big(\md{t+\nu\lambda},\cdot\big)\big\vert\J_t\big].$
If $P$ is ergodic with respect to $\theta_{\kappa_q(t)}$ then limit simplifies to
  $1/q\sum_{\nu=0}^{q-1}E\big[F\big(\md{t+\nu\lambda},\cdot\big)\big].$
}\end{ex}

\begin{ex}\label{ShiftsRet} Shifts as fibre transformations. {\rm
Consider the 2-parameter group define above \eqref{spc:twoShifts}.

\noindent
(i) If $\t$ is ergodic then, for any function $F\in\L^1(\ov \W,\ov \F,\ov P),$
\begin{eqnarray*}
  \frac{1}{n}\sum_{i=0}^{n-1}F\big(\t^i(t),
  \vartheta_{\kappa_i^{(1)}(t)\,v_1+\kappa_i^{(2)}(t)\,v_2}\big)
  \convinfty{n}\ov E[F]
   \qquad \ov P\text{-almost surely and in } \L^1(\ov P).
\end{eqnarray*}

\noindent
(ii) If $\t$ is periodic with period $q$ then,  for any function $F\in\L^1(\ov \W,\ov \F,\ov P),$
\begin{eqnarray*}
  \frac{1}{n}\sum_{i=0}^{n-1}F\big(\t^i(t),
  \vartheta_{\kappa_i^{(1)}(t)\,v_1+\kappa_i^{(2)}(t)\,v_2}\w\big)
  \convinfty{n}
  \frac{1}{q}\sum_{\nu=0}^{q-1}E\big[F\big(\t^\nu(t),
 \vartheta_{\kappa_\nu^{(1)}(t)\,v_1+\kappa_\nu^{(2)}(t)\,v_2}(\cdot)\big)\big\vert\, \J_t\big],
\end{eqnarray*}
for $\mu$-almost all
$t\in M,$ and for $P$-almost all $\w\in\W$ and in $\L^1(P).$
If  $\vartheta_{\kappa_q^{(1)}(t)\,v_1+\kappa_q^{(2)}(t)\,v_2}$ is ergodic with respect to $P,$
for $\mu$-almost all $t\in M,$ then the limit simplifies to
$  1/q\sum_{\nu=0}^{q-1} E\big[F\big(\t^\nu(t),\,\cdot\,\big)\big].$
}\end{ex}

\section{Uniform convergence}\label{unif}

This section addresses the question of {\it sure} convergence
with respect to the first parameter.
In addition to the assumptions at the beginning of Section \ref{spc}
we suppose that $M$ is a compact separable metric space endowed
with metric $d,$ and $\B$ is the Borel $\sigma$-algebra on $M$
for the topology induced by $d.$
Recall that the convergence of ergodic averages need not be
true {\it everywhere},
even if we are in a compact topological space and both the transformation
and the function are continuous. 
Which conditions guarantee sure convergence in the
first parameter? We will be asking
a little more than this, namely about {\it uniform} convergence in $t.$
We are interested in results of the type
\begin{equation}\label{unif:conv}
  \frac{1}{n}\sum_{i=0}^{n-1}
  F\big(\t^i(t),\theta_{\kappa_i(t)}\w\big)\convinfty{n}
  \ov E[F\vert\J](t,\w)\qquad\text{uniformly in }t\in M\text{ and in }\L^1(P).
\end{equation}
We further investigate when (\ref{unif:conv})
takes place $P$-almost surely, i.e., for $P$-almost all $\w\in\W,$
\begin{eqnarray}\label{UnifConvP}
 \lim_{n\to\infty}\;\sup_{t\in M}\bigg\vert
 \frac{1}{n}\sum_{i=0}^{n-1}F\big(\t^i(t),\theta_{\kappa_i(t)}\w\big)
 - \ov E[F\vert\J](t,\w)\,\bigg\vert=0.
\end{eqnarray}

Again, we consider two different cases: when $\t$ is ergodic
and when it is periodic.  In the second case we proceed as in the proof 
of Corollary \ref{Retper} and obtain the following uniform version.

\begin{cor}\label{Retperunif}
Let $\t$ be periodic with $q\in\N.$ Assume
$F\in\L^1(\ov\W,\ov\F,\ov P)$ with $F(t,\cdot)\in\L^1(\W,\F,P)$
for all $t\in M.$  Then, for $P$-almost all $\w$ and in $\L^1(P),$
\begin{equation*}
  \frac{1}{n}\sum_{i=0}^{n-1}F\big(\t^i(t),
  \theta_{\kappa_i(t)}\,\cdot\big)\convinfty{n}
            \frac{1}{q}\sum_{\nu=0}^{q-1}E\big[F\big(\t^\nu(t),
  \theta_{\nu_q(t)}\,\cdot\big)\big\vert\J_t\big]
  \qquad\text{uniformly in }t\in M,
\end{equation*}
where $\J_t$ denotes
the $\sigma$-algebra of $\theta_{\kappa_q(t)}$-invariant sets in $\F.$
If $P$ is ergodic with respect to $\theta_{\kappa_q(t)},$ for all $t\in M,$ the limit equals 
$1/q\sum_{\nu=0}^{q-1} E\big[F\big(\t^\nu(t),\,\cdot\,\big)\big].$
\end{cor}

The ergodic case is more delicate.  We begin with careful
investigations of ergodic theorems on the single spaces,
which will later be combined to derive a result on the product space.  
In the base we are dealing with a transformation on a compact and metrizable space. 
We now recall and refine some of the existing results about uniform convergence
in this situation. 
Motivated by applications in information theory.
we put particular emphasis on extending uniform convergence 
results to the class of Riemann-integrable functions.
An example for a function that is Riemann-integrable but not
continuous occurs in the proof of a directional Shannon-MacMillan theorem
for random fields (cf.~\cite{Bre01SMcM}).

The classical example for an ergodic theorem that gives
a statement about uniform convergence is the one by Weyl.  In its simplest
form, it says that the averages of a continuous function along
the orbit of an irrational translation on the circle converge uniformly
to the integral of the function.
To prove Weyl's theorem, Krengel (cf.~Theorem 2.6 in
Paragraph 1.2.3 in \cite{Kre85}) uses an Arzela-Ascoli
technique which we will make use of at the end of this section.
Under the assumption that $\t:M\rightarrow M$ is continuous
and that $f$ is a function on $M,$ such that the functions
$ F_n:=1/n\sum_{i=1}^{n-1} F\circ\t^i\ (n\in\N)$
are equicontinuous on $M,$ Krengel's theorem states that the
convergence in Birkhoff's ergodic theorem is uniform in $t.$
Together with the following Lemma, it yields Weyl's theorem.

\begin{lemma}\label{EquBed}
Let $f$ be a continous function on $M,$ and
$\t:M\rightarrow M$ Lipschitz-continuous with Lipschitz
constant $c\le 1.$ Then the functions $F_n\ (n\in\N)$ defined above are equicontinuous.
\end{lemma}

\pf We have to show that for
every $\varepsilon>0$ there is a $\delta>0$ such that for all
$n\in\N$ and all $s,t\in M$ with $d(s,t)<\delta,$
  $1/n\big\vert\sum_{i=0}^{n-1}f(\t^i(s))-f(\t^i(t))\big\vert
  <\varepsilon.$
Fix $\varepsilon>0.$ We will show
that there is a $\delta>0,$ such that
$\vert f(\t^i(s))-f(\t^i(t))\big\vert <\varepsilon$
for all $d(s,t)<\delta.$ Since $M$ is compact, $f$ must be
uniformly continuous, i.e., there is a $\delta>0$ such that
for all $x,y\in M$ with $d(x,y)<\delta$ we have
$|f(x)-f(y)|<\varepsilon.$
By assumption,
$d(\t(s),\t(t))\le c\,d(s,t)$ for all $s,t\in M,$
and therefore, $d(\t^i(s),\t^i(t))\le c^i\,d(s,t)\le d(s,t)$
for all $s,t\in M,$ and for all $i\in\N.$
\qed\pfskip

We shall ask whether we could replace the assumption of continuity
of the function $f$ in Weyl's theorem by a weaker condition. It is
certainly not true for all measurable functions, which can be seen
in a simple example: Fix $t_0\in\T.$ Its orbit under $\t$ is the
set ${\cal O}:=\{\t^n(t_0)\vert n\in\N_0\,\}.$ 
For the function $f:=1_{\cal O},$ 
the ergodic averages converge to 1, for all $t\in {\cal O,}$
but $\int_0^1 1_{\cal O}(t)\,dt=0.$  

This question of uniform convergence has a connection with unique
ergodicity (cf.~e.g., Chapter 4.1.e.\ of \cite{HaK95} or 
Theorem 6.19 in \cite{Wal82}).  A continuous transformation $\t$ of a 
compact metrizable space is called {\rm uniquely ergodic} if it has only 
one invariant Borel measure.
It can be shown that this measure must be ergodic, which implies
that the ergodic averages of an integrable function converge
almost surely to a constant. 
The Lebesgue measure is the only probability measure on $(\T,\B),$
which is invariant with respect to rotations of the circle,
so it must be uniquely ergodic.
Oxtoby extended Weyl's theorem to the situation where $\t$ is a
uniquely ergodic transformation on a compact metric space.
It states uniform convergence of the ergodic  averages of a
continuous function.  
Note that, conversely, uniform convergence does
not imply the continuity of the function.
(Further conditions for this would be needed, such as 
topological transitivity of $\t$ or constancy of the limit.)

Below Theorem 2.7 in Chapter 1 of \cite{Kre85}, Krengel mentions
that Weyl's theorem is sometimes spelled out for
to the class of Riemann-integrable functions.
Actually, it was proved by de Bruijn and Post
\cite{BrP68} that the function is
Riemann-integrable if and only if the convergence is uniform. 
This also follows from our next proposition. We ask the following question: 
Considering uniform convergence of the ergodic
averages along a continuous transformation on a compact
real interval, can we pass automatically from the class of
continuous functions to the class of functions which are
integrable in the sense of Riemann?

\begin{prop}\label{WeylR}
Let $a,b\in\R,$ $a<b,$ $\mu$ a measure on $([a,b],\B([a,b]))$
which is absolutely continuous with respect to Lebesgue measure,
with a continuous density. Let $\t:[a,b]\rightarrow [a,b]$ be
continuous. Assume that for any continuous function $f$ on $[a,b]$
\begin{eqnarray}\label{WeylR1}
  \ErgAv{f}{\t}\convinfty{n}\int_a^b f\,d\mu\qquad
  \text{uniformly.}
\end{eqnarray}
Then the convergence holds as well for any function which
is integrable in the sense of Riemann.
\end{prop}

\noindent The proof is carried out using a common sandwich argument
(cf.~e.g.,\ in Chapter 4.1.e.\ of \cite{HaK95}).  Applying the proposition
to the situation of Oxtoby's Theorem yields

\begin{cor}\label{UnifR} Let $a,b\in\R,$
$a<b,$ and $\t:[a,b]\rightarrow [a,b]$ continuous and uniquely
ergodic with invariant measure $\mu.$ Assume that $\mu$ is
absolutely continuous with respect to Lebesgue measure, with a
continuous density. Then, for any function $f$ on $M$ which is
integrable in the sense of Riemann,
\begin{displaymath}
  \ErgAv{f}{\t}\convinfty{n}\int_0^1 f\,d\mu\quad\quad\text{uniformly.}
\end{displaymath}
\end{cor}

Now we focus on studying the convergence in the fibre.
Fix $t\in M,$ and define a function on $\W$ by $f(\w):=F(t,\w).$ 
This reduces the ergodic averages of the skew product to
$1/n\sum_{i=0}^{n-1}f\big(\theta_{\kappa_i(t)}\w\big),$
which we can view as a sort of
ergodic average along the subsequence $(\kappa_i(t))_{i\in\N}.$
Recalling a classical result about  $\L^p$-convergence of ergodic
averages along subsequences, there is the following characterization.

\begin{thm}\label{BHthm}{\rm (Blum \& Hanson)}
Let $T$ be a transformation on $(\W,\F).$ Suppose that $T$ is
invertible and that both, $T$ and $T^{-1}$ preserve $P.$ Then $P$
is strongly mixing with respect to $T$ if and only if for all $p,\
1\le p<\infty,$ every strictly increasing sequence $(m_i)_{i\in
\N}$ of integers, and every function $f\in\L^p(\W,\F,P),$
\begin{equation*}
  \frac{1}{n}\sum_{i=0}^{n-1}f\circ T^{m_i}
  \convinfty{n} E[f]\quad\quad\text{in }\L^p(P).
\end{equation*}
\end{thm}
\noindent The key to the proof of Blum and Hanson's theorem is the following

\begin{lemma}\label{BHlem}
Under the assumptions of Theorem \ref{BHthm} and supposing that
$P$ is strongly mixing with respect to $T$ we have for all $A\in\F,$
for every strictly increasing sequence $(k_i)_{i\in\N},$
\begin{equation*}
  \frac{1}{n}\sum_{i=0}^{n-1}1_A\circ T^{k_i}\convinfty{n}
   P(A)\quad\quad\text{ in }\L^2(P).
\end{equation*}
\end{lemma}

Our next step is to carry over Blum and Hanson's theorem to the
case of a $d$-parameter group of transformations $(\theta_k)_{k\in\Z^d},$
at the place of iterates of one transformation $T.$
What we need is a condition on the $\Z^d$-valued 
sequence $(k_i)_{i\in \N}$ which replaces the strict
monotonicity imposed on $(m_i)_{i\in\N}.$ With an eye toward later
applications on the product space we generalize the result further
by showing that the $\L^2$-convergence takes place uniformly over a family of functions,
indexed by a set $I.$
Recall that $\|\cdot\|$ denotes the maximum norm in $\Z^d.$

\begin{lemma}\label{BHtyplem}
Assume that $P$ is strongly mixing with respect to
$(\theta_k)_{k\in\Z^d},$ and let $(k_n(t))_{n\in\N}$ $(t\in I)$ be a
family of sequences with values in $\Z^d$ that fulfill, for all $m\in\N,$
\begin{equation}\label{BHA}
  \lim_{n\to\infty}
  \frac{1}{n^2}\sup_{t\in I}\,\big|\big\{1\le i,j\le n\,
  \big\vert\,\|k_i(t)-k_j(t)\|\le m\big\}\big|=0.
\end{equation}
Then for all $A\in\F,$
\begin{equation*}
 \sup_{t\in I}\bigg\|
 \frac{1}{n}\sum_{i=0}^{n-1}1_A\circ\theta_{k_i(t)}-
 P(A)\bigg\|_{\L^2(P)}^2 \convinfty{n}0.
\end{equation*}
\end{lemma}

\pf
For every $A\in\F$ and $t\in I$ we obtain by simple calculations,
\begin{align}
   \bigg\|\frac{1}{n} & \sum_{i=0}^{n-1}  1_A\circ\theta_{k_i(t)}-P(A)\,
   \bigg\|_{\L^2(P)}^2
     =\int_{\W}\frac{1}{n^2}\sum_{i,j=0}^{n-1}(1_A\circ\theta_{k_i(t)}-P(A))
    (1_A\circ\theta_{k_j(t)}-P(A))\,dP\nn\\
      &=\frac{1}{n^2}\sum_{i,j=0}^{n-1}\bigg[\int_{\W}\big(
   1_A\circ\theta_{k_i(t)}\, 1_A\circ\theta_{k_j(t)}\big)\,dP - P(A)
   \int_{\W}\big(1_A\circ\theta_{k_i(t)}
   + 1_A\circ\theta_{k_j(t)}\big)\,dP + P(A)^2\bigg]\nn\\
   \label{BH1}
 &=\frac{1}{n^2}\sum_{i,j=0}^{n-1}\Big(P\big(\theta_{k_i(t)}^{-1}A
   \cap\theta_{k_j(t)}^{-1}A\big)-P(A)^2\Big)
   \le\frac{1}{n^2}\sum_{i,j=0}^{n-1}\Big|
    P\big(\theta_{k_i(t)-k_j(t)}^{-1}A\cap A\big) - P(A)^2\Big|.
\end{align}
Fix $\varepsilon>0.$ Due to the mixing condition there is
an $m\in\N$ such that
\begin{equation*}
  \Big| P\big(\theta_{k_i(t)-k_j(t)}^{-1}A\cap A\big)
  - P(A)^2\Big|
  <\frac{\varepsilon}{2}\quad\quad\text{for all }i, j\in\N
  \text{ with }\|\,k_i(t)-k_j(t)\|> m,
\end{equation*}
and
\begin{equation*}
  \Big| P\big(\theta_k^{-1}A\cap A\big) - P(A)^2\Big|
  <\frac{\varepsilon}{2}\quad\quad\text{for all }k\in\Z^2
  \text{ with }\|k\|\ge m.
\end{equation*}
By assumption (\ref{BHA}) there is a $n_0\in\N$ such that
\begin{equation*}
  \frac{1}{n^2}\sup_{t\in I}\,\big|\big\{1\le i,j\le
  n\,\big\vert\,\|k_i(t)-k_j(t)\|\le m\big\}\big|
  <\frac{\varepsilon}{2}
  \quad\text{ for all }n\ge n_0.
\end{equation*}
Applying the last two inequalities to (\ref{BH1})
yields for all $n\ge n_0$
\begin{align*}
  \sup_{t\in I}
& \bigg\|\frac{1}{n^2}\sum_{i=0}^{n-1}
  1_A\circ\theta_{k_i(t)} -P(A)  \,
  \bigg\|_{\L^2(\W,\F,P)}^2   < \varepsilon,
\end{align*}
\nopagebreak
and the assertion of the lemma follows by letting $\varepsilon$ to $0.$
\qed\pfskip

\begin{thm}\label{BHtk}
Assume that $P$ is strongly mixing with
respect to $(\theta_k)_{k\in\Z^d}.$ Let $(k_n(t))_{n\in\N}$ $(t\in I)$ be a
family of sequences with values in $\Z^d,$ for which for all $m\in\N,$
\begin{equation*}
  \lim_{n\to\infty}\frac{1}{n^2}\,
  \sup_{t\in I}\big|\big\{1\le i,j\le
   n\,\big\vert\,\|k_i(t)-k_j(t)\|\le m\big\}\big|=0.
\end{equation*}
Then for $1\le p<\infty$ and for any $f\in\L^p(\W,\F,P),$
\begin{equation}\label{BHtkSta}
 \frac{1}{n}\sum_{i=0}^{n-1}f\circ\theta_{k_i(t)}\convinfty{n}E[f]
 \quad\text{ in }\L^p(P), \text{  uniformly in }t\in I.
\end{equation}
\end{thm}

\pf As an immediate consequence of the preceeding lemma,  
\eqref{BHtkSta} is true for $p=2$ for any simple function $g$ on $(\W,\F).$ 
By a standard argument (cf., e.g., Lemma 4 in \cite{BlH60}),
this convergence holds as well in $\L^p(P),$ for
$1\le p<\infty.$ Finally, for any function
$f$ in $\L^p(P),$ decomposition into positive and negative parts,
$\L^p(P)$-approximation by simple functions, and monotone
convergence yields (\ref{BHtkSta}).
\qed\pfskip

We will now combine the approaches developed separately for the base
and the fibre transformation to derive a result about uniform convergence
in the base and $\L^p$-convergence in the fibre.  A crucial ingredient is a condition 
that regulates the effects of the coupling sequence $(\kappa_n(t))_{n\in\N}.$ 

\begin{thm}\label{UnifL}
Let $\t:M\rightarrow M$ be continuous and uniquely ergodic, and
suppose that $P$ is strongly mixing
with respect to the group of transformations $(\theta_k)_{k\in\Z^d}.$ Let
$\kappa: M\rightarrow \Z^d$ be $\B$-measurable such that,
\begin{equation}\label{kappaCond}
  \lim_{n\to\infty}\frac{1}{n^2}\,
   \sup_{t\in M}\big|\big\{1\le i,j\le
   n\,\big\vert\,\|\kappa_i(t)-\kappa_j(t)\|\le
   m\big\}\big|=0
\end{equation}
for all $m\in\N.$ Let be $1\le p<\infty.$ Then for every $\F$-measurable function
$F$ on $\ov\W$ such that $\sup_{t\in M} \vert\, F(t,\cdot\,)\,\vert$ is in
$\L^{p}(\Omega,\F,P)$ and $F(\,\cdot\,,\omega)$ is continuous on $M$
for P-almost every $\omega$,
\begin{equation}\label{UnifLClaim}
  \frac{1}{n}\sum_{i=0}^{n-1}F\big(\t^i(t),
  \theta_{\kappa_i(t)}\,\cdot\big) \convinfty{n}
  \ov E[F]
  \qquad\mbox{in }\L^p(P)\mbox{  and uniformly in }t\in M.
 \end{equation}
\end{thm}

\pf
We first prove the theorem in the case when $F$ is the indicator of a
set of the form $U\times A$, where $U$ is the intersection of
finitely many metric balls in $M$ or their complements, and $A\in\F.$
By (\ref{ItS}), the expression
\begin{equation}\label{UnifLDiff}
 \bigg\|\frac{1}{n}\sum_{i=0}^{n-1}F\big(\t^i(t),
  \theta_{\kappa_i(t)}\,\cdot\big) -\ov E[F]\bigg\|_{\L^2(P)}
\end{equation}
then transforms to
\begin{equation}\label{BHprdlem1}
\begin{split}
   \frac{1}{n^2} & \sum_{i,j=0}^{n-1}\int_{\W}1_U(\t^i(t)) 1_U(\t^j(t))
   1_A(\theta_{\kappa_i(t)}\w)1_A(\theta_{\kappa_j(t)}\w)\,P(d\w) \\
    - & \frac{1}{n^2}\sum_{i,j=0}^{n-1}\mu(U)P(A)\bigg[  
    \int_{\W}1_U(\t^i(t))1_A(\theta_{\kappa_i(t)}\w)\,P(d\w)
    + \int_{\W}1_U(\t^j(t))
   1_A(\theta_{\kappa_j(t)}\w)\,P(d\w)\bigg] \\ 
    \vspace{-2mm}
    + &\  \mu(U)^2 P(A)^2.
\end{split}
\end{equation}
By $ P\big(\theta_{\kappa_i(t)}^{-1}A\cap\theta_{\kappa_j(t)}^{-1}A\big) =
  P\big(\theta_{\kappa_i(t)-\kappa_j(t)}^{-1}A\cap A\big),$
the first addend equals
\begin{equation*}
  \frac{1}{n^2}\sum_{i,j=0}^{n-1}1_U(\t^i(t))1_U(\t^j(t))
    P\big(\theta_{\kappa_i(t)-\kappa_j(t)}^{-1}A\cap A\big).
\end{equation*}
It may be replaced by
\begin{equation}\label{BHprdlem2}
  \frac{1}{n^2}\sum_{i,j=0}^{n-1}1_U(\t^i(t))1_U(\t^j(t)) P(A)^2
\end{equation}
without affecting the asymptotic behavior of (\ref{BHprdlem1}) as can be seen as follows.
We may bound
\begin{align*}
  \bigg\vert&\frac{1}{n^2}\sum_{i,j=0}^{n-1}1_U(\t^i(t))1_U(\t^j(t))
   \Big(P\big(\theta_{\kappa_i(t)-\kappa_j(t)}^{-1}A\cap A\big)
    -P(A)^2\Big)\bigg\vert\nn\\
  &\le  \frac{1}{n^2}\sum_{i,j=0}^{n-1}
    \Big\vert P\big(\theta_{\kappa_i(t)-\kappa_j(t)}^{-1}A\cap A\big)
    -P(A)^2\Big\vert.
\end{align*}
Now, we argue as in the second part of the proof of Lemma
\ref{BHtyplem}, replacing the sequence $(k_i)_{i\in\N}$ by
$(\kappa_i(t))_{i\in\N},$ and using assumption (\ref{kappaCond})
instead of (\ref{BHA}). This proves that the difference
created by the change \eqref{BHprdlem2}
converges to $0$ uniformly with respect to $t.$

Since the term in the rectangular brackets in the second addend in
(\ref{BHprdlem1}) equals
   $ 1_U(\t^i(t))P(A)+1_U(\t^j(t))P(A),$
the whole expression (\ref{BHprdlem1}) simplifies to
\begin{equation*}
   \frac{1}{n^2}\sum_{i,j=0}^{n-1}\Big(
    1_U(\t^i(t))1_U(\t^j(t))
    -\mu(U)\big(1_U(\t^i(t))+1_U(\t^j(t))\big)
    +\mu(U)^2\Big) P(A)^2,
\end{equation*}
which can be further reduced to
$ \bigg(1/n\sum_{i=0}^{n-1}1_U(\t^i(t))-\mu(U)\bigg)^2 P(A)^2.$
Since $\t$ is uniquely ergodic and $\mu(\partial U)=0,$
Corollary 4.1.14 in \cite{HaK95} tells us that the first factor converges
to $0$ uniformly in $t,$
which concludes the first part of the proof.
To pass from $\L^2$-convergence to general $\L^p,$ use again
a standard argument (for instance, Lemma 4 in \cite{BlH60}).

Now we let $F$ be a general function, satisfying the conditions of
the theorem.  We need to find for every positive
$\epsilon$ a sequence of metric balls $U_i$ in $M$ and $A_i\in\F$, with
real numbers $a_{i}$ such that, for all $t\in M$,
$\| F(t,\cdot)- I(t,\cdot) \|_{p}< \eps,$ where
$I=\sum_{i=1}^{n} a_{i} 1_{U_{i}\times A_{i}}.$
It will then follow that
\begin{align*}
\bigg\|\frac{1}{n} \sum_{i=0}^{n-1}
F\bigl(\tau_{i}(t),\theta_{\kappa_{i}(t)} \,\cdot\, \bigr) -
\frac{1}{n} \sum_{i=0}^{n-1}
I\bigl(\tau_{i}(t),\theta_{\kappa_{i}(t)} \,\cdot\, \bigr) \bigg\|_{p}
 \leq \frac{1}{n} \sum_{i=0}^{n-1}
  \bigl\|F\bigl(\tau_{i}(t),\,\cdot\, \bigr) -
  I\bigl(\tau_{i}(t), \,\cdot\, \bigr) \bigr\|_{p}
  < \eps.
\end{align*}
For $\omega\in\Omega$ and $c>0$, let $\delta(c,\omega)$ be the modulus
of continuity for the function $F(\,\cdot\,,\omega)$.
Define the sets
$$
M_k=\sup_{t\in M}\{\w\,\vert\, |F(t,\omega)|\le k\}
\quad\text{and}\quad
D_k(\eps)=\{\w\,\vert\,\delta(1/k,\omega)\le\eps\}
\quad\ (k\in\N).
$$
Then the sequence of functions on $\Omega,$
$
F_{k}(\omega):=\sup_{t\in M}|F|^{p}(t,\omega)\,
1_{D_k(\eps/4)^c\cup M_k^c}(\w)\ (k\in\N),
$
is bounded by $\sup_{t\in M}|F|^{p}(t,\omega)$, which is integrable, and converges to
$0$ for every $\omega$.  By the bounded
convergence theorem, the integral of $F_{k}$ converges to 0 as $k$
goes to infinity.  Choose a $k$ such that $\int F_{k}\,P(d\w)\leq (\eps/2)^{p}$.
Since $M$ is compact, we may find a finite sequence
$t_{1},\dots,t_{r}\in M$ such that the balls of radius $1/k$ around
these centers cover $M$.  We also define a sequence of real numbers
$-k-1=s_{0}<\cdots<s_{r'}=k$ such that the difference between any two
successive elements is less than $\eps/8$.  Now we define a collection
of sets $U_{i,j}$ and $A_{i,j}$
indexed by $r\times r'$.  We start with $U_{i,j}$ as the  ball of radius
$1/k$ around $t_{i}$, and then remove the intersections, so that the
$U_{i,j}$ is the same for all $j$, and running through $1\leq i\leq
r$ yields a disjoint cover of $M$.  The sets $A_{i,j}$ are defined by
$$
A_{i,j}=\big\{ \w\,\big\vert\, s_{j-1}<F(t_{i},\omega)\leq s_{j}\big\}
\cap D_k(\eps/8) \cap M_k.
$$
Let $a_{i,j}=s_{j}$. We throw in one additional product set,
$U_{0}=M$ and
$A_{0}=D_k(\eps/8)^c\cup M_k^c$
with $a_{0}=0$, and define the simple function $I(t,\omega)$ as
indicated above.
Then for any $t\in M,$
\begin{align*}
    \Bigl  \| F(t, & \cdot) -I(t,\cdot)\Bigr\|_{p}  
    =
     \bigg(
      \int \big\vert  F(t,\omega)-I(t,\omega) \big\vert^{p} 
      1_{D_k(\eps/8)}\,1_{M_k}\,P(d\omega)
    \bigg)^{1/p} 
    + \: \bigg(\int F_{k}(\omega)\, P(d\omega)
          \bigg)^{1/p}.
\end{align*}
We already assumed (in defining $k$) that the second term is smaller
than $\eps/2$.  For every $t$, there is a unique pair $i,j$ such that
$t\in U_{i,j}$ and $\omega\in A_{i,j}$.  By construction,
$I(t,\omega)=s_{j}$, so the integrand in the first term is bounded by
$
2^{p}\bigl|F(t_{i},\omega)-F(t,\omega) \bigr|^{p} +
  2^{p}\bigl|F(t_{i},\omega)-s_{j} \bigr|^{p} .
  $
This in turn is bounded by $2^p(\eps/4)^{p}< (\eps/2)^{p}$, since
$\omega$ is not in $A_{0}$ and $d(t_{i},t)<1/k$, completing the proof.
\qed\pfskip

Proposition \ref{WeylR} yields the following version for Riemann-integrable functions.

\begin{cor}\label{UniLRMak}
Let $a,b\in\R,$ $a<b,$ and $\t:[a,b]\rightarrow [a,b]$ be continuous
and uniquely ergodic with invariant measure $\mu,$ and
assume that $\mu$ is absolutely continuous
with respect to Lebesgue measure, with continuous density.
For $P$ and $\kappa$ assume the same as in Theorem \ref{UnifL}.
Let $F\in L^p([a,b]\times\W,\B\otimes\F,\ov P)$
be Riemann-integrable with respect to the first
variable. Then we have
\begin{eqnarray*}
 \lim_{n\to\infty}\;\sup_{t\in M}\bigg\|
 \frac{1}{n}\sum_{i=0}^{n-1}
 F\big(\t^i(t),\theta_{\kappa_i(t)}\cdot\big)
 - \ov E[F]\,\bigg\|_{\L^p(P)}=0.
\end{eqnarray*}
\end{cor}

Our next goal is to derive a statement about $P$-almost sure convergence
rather than $\L^1(P)$-convergence in the fibre.  Further conditions
on $F$ are needed.  $P$-almost sure convergence of ergodic theorems
involving weights or subsequences is a very subtle question (cf., e.g. \cite{BeL85}).
Choosing a function $F$ which is 
constant in $\w$ and considering Lemma \ref{EquBed}
suggests that we need an equicontinuity assumption in $t.$ 
Note also the additional assumptions that non-empty open sets on $M$
have positive mass under $\mu.$

\begin{thm}\label{UnifP}
Let $\mu$ be a $\t$-invariant measure on $(M,\B),$
such that any non-empty open subset $U$ of $M$ has $\mu(U)>0.$
Let $\t:M\rightarrow M$ be continuous and $F$ a function
on $M\times\W,$  for which $F(t,\cdot\,)\in\L^1(P)$
for all $t\in M,$  and the sequence of functions
$$
\Bigg(\frac{1}{n}\sum_{i=0}^{n-1}F(\t^i(\cdot\,),
\theta_{\kappa_i(\cdot\,)}\w)\Bigg)_{n\in\N}
$$
is equicontinuous on $M,$ for all $\w\in\W.$
Then, for $P$-almost all $\w\in\W,$
\begin{eqnarray}\label{UnifPBeh}
  \sup_{t\in M}\bigg\vert
 \frac{1}{n}\sum_{i=0}^{n-1}F\big(\t^i(t),\theta_{\kappa_i(t)}\w\big)
 -\ov E[F\vert\J](t,\w)\,\bigg\vert\convinfty{n}0.
\end{eqnarray}
\end{thm}
\pf
We may assume without loss of generality that $E[F\vert\J]=0.$
The general case can be reduced to this by subtracting $E[F\vert\J]$
on both sides and making use of the invariance of $E[F\vert\J]$ under $S.$
The first step is to construct a countable dense set $M_1\subset M$
and a set $N_1\subset\W$ with $P(N_1)=0$ such that
\begin{eqnarray}\label{UnifP1}
  \frac{1}{n}\sum_{i=0}^{n-1}F\circ S^i(t,\w)\longrightarrow 0
  \quad\quad\text{ for all $t\in M_1$ and all }\w\in\W\setminus N_1.
\end{eqnarray}
Since $M$ is compact, the conditions on $F$ assure that
$F\in\L^1(\ov\W,\ov\F,\ov P),$ and therefore by (\ref{BirkS})
there is a set $M_1\subset M$ with $\mu(M_1)=1$ such that for
any $t\in M_1$ there is a set $N(t)\subset\W$ with $P(N(t))=0$ and
\eqref{UnifP1} holds for all $\w\in\W\setminus N(t).$
${\widetilde M}_1$ is dense in $M$ because its complement has measure zero
with respect to $\mu$ and therefore, by assumption,
contains no non-empty open subsets.
Since $M$ is separable we can find a countable dense subset $C\subset M,$
and because ${\widetilde M}_1$ is dense in $M,$ we can approximate any
$x\in C$ by a sequence $(a_j(x))_{j\in\N}$ with $a_j(x)\in {\widetilde M}_1$
for all $j\in\N.$ 
$M_1:=\bigcup_{x\in C}\bigcup_{j\in\N}a_j(x)$ defines a countable dense subset
of $M,$ and $N_1:=\bigcup_{t\in M_1} N(t)$ defines a subset of
$\W,$ which fulfills (\ref{UnifP1}).  This completes the first step.

For the next step, choose $s\in M$ and
fix $\varepsilon>0.$ By equicontinuity, there is a set
$N_0\subset\W$ with $P(N_0)=0$ and a
$\delta>0$ such that for all $r,t\in M$ with $d(r,t)<\delta$
\begin{eqnarray}\label{UnifP3}
  \bigg\vert\frac{1}{n}\sum_{i=0}^{n-1}F\circ S^i(r,\w)-F\circ
  S^i(t,\w)\bigg\vert
  <\frac{\varepsilon}{2}\quad\text{for all }n\in\N\text{ and all }
  \w\in\W\setminus N_0.
\end{eqnarray}
Define $N:=N_0\cup N_1$ and fix $\w\in\W\setminus N.$
Since $M_1$ is dense in $M$ we can find a $t\in M_1$ with $d(s,t)<\delta,$
and by (\ref{UnifP1}) there is an $n_1\in\N$ such that
\begin{eqnarray*}
  \bigg\vert\ErgAv{F}{S}(t,\w)\bigg\vert <\frac{\varepsilon}{2}
  \quad\quad\text{ for all }n\ge n_1.
\end{eqnarray*}
Combining the last two inequalities leads the desired
\begin{eqnarray*}
  \bigg\vert\ErgAv{F}{S}(s,\w)\bigg\vert <\varepsilon
  \quad\quad\text{ for all }n\ge n_1\text{ and all }
  \w\in\W\setminus N.
\end{eqnarray*}

For uniform convergence w.r.t.~the first variable,
we use a standard compactness argument.  $M$ 
can be covered by a finite number $m$ of $\delta$-neighborhoods
in $M$, which centers are denoted by $s_1,...,s_m.$ Applying the
reasoning of the last step to each of the $s_1,...,s_m$ we 
find $n_0\in \N$ such that
\begin{eqnarray*}
  \bigg\vert\ErgAv{F}{S}(s_k,\w)\bigg\vert <\frac{\varepsilon}{2}
  \quad\quad\text{ for all }n\ge n_0, k\in\{1,...,K\},
  \text{ and }\w\in\W\setminus N.
\end{eqnarray*}
For an arbitrary $s\in M$ there exists $k\in\{1,...,K\}$ such that
$d(s,s_k)<\delta,$ and by (\ref{UnifP3}) we obtain\begin{eqnarray*}
  \frac{1}{n}\bigg\vert\sum_{i=0}^{n-1}F\circ S^i(s,\w)-F\circ
  S^i(s_k,\w)\bigg\vert
  <\frac{\varepsilon}{2}\qquad\text{  for all }n\in\N\text{ and all }
  \w\in\W\setminus N.
\end{eqnarray*}
Finally, the convergence (\ref{UnifPBeh}) follows by the last two inequalities.
\qed

We conclude the paper with an application of our results to the question that originally 
motivated this work.  Consider the ergodic averages of a function of a random field restricted 
to the staircase pattern defined by the approximation of a line by a lattice as in \eqref{Intro:ex.sq}.  
We actually the function to be multivariate.  In this context, this means that it may depend on
a finite number of such steps.
What can be said about the asymptotic behavior of the ergodic averages of such a function?

Let $P$ be a two-dimensional random field, that is, a probability measure on 
$\Omega=\Upsilon^{\Z^2}$ invariant w.r.t.~the group of shift transformations 
$(\vartheta_v)_{v\in\Z^2}$ (see just above Corollary \ref{ShiftsSProp} for details).
For a set $U\subseteq \Z^2$ let $\pi_U(\w):=\w(U)$ and $P_U:=P\circ \pi_U^{-1}$ defines a
probability distribution on $\Upsilon^{\vert U\vert}.$
Use $[x]$ denote the integer part of a real number $x,$ respectively.
Recall that $L_{\lambda,t}(z)=(z,[\lambda z+t])\ (z\in \Z)$ is an approximation of the line with slope 
$\lambda$ and $y$-intercept $t$ by elements of the lattice $\Z^{2}.$ For $z_1,z_2\in\Z,$
$L_{\lambda,t}(z_1,...,z_2)$ are the $z_1${\it th} to the $z_2${\it th} step.
We will use the short form $P_{\lambda,t,m}:=P_{L_{\lambda,t}(0,...,m-1)}.$
 
\begin{cor}\label{Unif:ex.staircase}   
Let $m\in\N$ and $\lambda\in [0,1].$
Let $f$ be a function on $\Upsilon^m$ which is integrable 
w.r.t.~$P_{\lambda,t,m}$ for all $t\in[0,1].$  

\noindent (i) 
Let $\lambda$ be rational and represented as
$\lambda=p/q$ for $p\in\Z$ and $q\in\N$ with no common divisor.
If $P$ is ergodic w.r.t.~the shift transformations then
\begin{equation*}
   \frac{1}{n}\sum_{i=0}^{n-1}f \big( \w (L_{\lambda,t}(i,...,i+m-1)) \big)
   \convinfty{n} 
   \frac{1}{q}\sum_{\nu=0}^{q-1} 
   \int_{\Upsilon^m} f(y)\,P_{\lambda,\t^\nu(t),m}(dy)
\end{equation*}
in $\L^1(P),$ for $P$-almost all $\w\in\Omega,$ and uniformly in $t\in [0,1].$

\noindent (ii) 
Let $\lambda$ be irrational. Assume further that
$f \big( \w (L_{\lambda,t}(0,...,m-1))\big)$ is Riemann-integrable with respect to $t\in [0,1],$
for $P$-almost all $\w\in\Omega.$
If $P$ is strongly mixing w.r.t.~the shift transformations then
\begin{equation*}
   \frac{1}{n}\sum_{i=0}^{n-1}f \big( \w (L_{\lambda,t}(i,...,i+m-1)) \big)
   \convinfty{n}  
   \int_0^1 \int_{\Upsilon^m} f(y)\,P_{\lambda,t,m}(dy)\,dt
\end{equation*}
in $\L^1(P)$ and uniformly in $t\in [0,1].$

\end{cor}

\pf 
Let $\t_{\lambda}$ be the rotation on the circle defined in Example \ref{TorTransl}, 
and $S_\lambda$ the skew product defined in \eqref{Intr:def.skew}.  
Use $\kappa(t):=(1,[t+\lambda])$ and $\kappa_n$ as defined in \eqref{ItS}.
It is easy to show that 
$\kappa_n(t)=L_{\lambda,t}(n)$ for all $n\in\N:$
For $n=1,$ it follows immediately from plugging in the definition,
and it remains to induce the statement from $n$ to $n+1.$
It is obvious for the first coordinate.  The second coordinate of $\kappa$
can be written as 
$\kappa_{n+1}^{(2)}(t)=\kappa_n^{(2)}(t)+\kappa\circ\tau_\lambda^n(t)
=[n \lambda + t] + [\tau_\lambda^n(t) + \lambda].$
The claim now follows from
$ [\t_{\lambda}^n(t)+\lambda]
= [t+n\lambda - [t+n\lambda] + \lambda]
=-[t+n\lambda] + [t+(n+1)\lambda].$

This implies that the iterates of the skew product are of the form
$S_{\lambda}^n(t,\omega)=\big(\t_\lambda^n(t),\theta_{L_{\lambda,t}(n)}\omega\big)$
$(n\in\N)$ capturing the lattice approximation of the line.
Using the second equality in \eqref{ItS}, it also follows that
$L_{\lambda,t}(n+u)=L_{\lambda,t}(n)+L_{\lambda,\t^n(t)}(u)$ for all $n,u\in\N_0.$
We thus get
$ L_{\lambda,t}(i,...,i+m-1) = L_{\lambda,t}(i) + L_{\lambda,\t^i(t)}(0,...,m-1).$
Introducing the function $F_\lambda(t,\w):=f\big(\w,L_{\lambda,t}(0,...,m-1)\big)$ we obtain
\begin{align}\label{Ex:RewrittenAs}
   \frac{1}{n}\sum_{i=0}^{n-1} & f \big( \w (L_{\lambda,t}(i,...,i+m-1)) \big) \nonumber \\
& = \frac{1}{n}\sum_{i=0}^{n-1}f \circ \vartheta_{L_{\lambda,t}(i)} 
   \big( \w (L_{\lambda,\t^i(t)}(0,...,m-1)) \big)
= \frac{1}{n}\sum_{i=0}^{n-1} F_\lambda\circ S_{\lambda}(t,\w).
\end{align}

Consider case (i).  With the representation $\lambda=p/q$ as above, $\t_\lambda$ is periodic with $q.$  
By Corollary \ref{Retperunif} and \eqref{Ex:RewrittenAs}, the averages on the left side in (i)
converge uniformly in $t,$ for $P$-almost all $\w$ and in $\L^1(P).$  
Since $P$ is ergodic with respect to $\theta_{\kappa_q(t)}$ for all $t\in M,$ the limit equals 
$1/q \sum_{\nu=0}^{q-1} \int_\W F_\lambda\big(\t_\lambda^\nu(t),\,\cdot\,\big)\,P(d\w)$
which simplifies to the expression on the right hand side of(i).

Consider case (ii).  For an irrational $\lambda,$ $\t_{{\lambda}}$ is uniquely ergodic.
Since $\|\kappa_n(t)\|$ (with $\|\,\cdot\,\|$ for the maximum norm) is bounded
from below by $n,$ the sequence tends to infinity as $n$ goes to infinity.  
Corollary \ref{ShiftsSProp} with $v_1=(1,0)$ and $v_2=(0,1)$
implies the ergodicity of $S_\lambda.$ 
It remains to verify condition \eqref{kappaCond}.
The latter easily follows from $\|\kappa_i(t)-\kappa_j(t)\|\ge | i-j |,$ and since
${1}/{n^2}\big\vert\big\{1\le i,j\le n\,\big\vert\,\vert i-j\vert\le m\big\}\big\vert$
converges to $0$ for all $m\in\N.$ 
Now, Corollary \ref{UniLRMak} applied to \eqref{Ex:RewrittenAs} implies that 
the averages in (ii) converge uniformly in $t,$ and in $\L^1(P).$  
The limit equals $\int_0^1 \int_\W F_\lambda(t,\,\cdot\,)\,P(d\w)\,dt$ which
simplify to the expression on the right hand side of (ii).
\qed

Note that $P$-almost everywhere convergence in (ii) can be derived as well, but 
requires additional conditions of the form stated in Theorem \ref{UnifP}.
For $m=1$ the limit actually simplifies to the integral of $f$ with respect to
the marginal distribution of $P$ in the origin.
In particular, it is independent of $t$, and it is the same in (i) and (ii).
The simplest interesting case is $m=2.$  
Let $P_{\mbox{flat}}:=P\circ\pi_{\{(0,0),(1,0)\}}^{-1}$ denote the marginal distribution
of $P$ on the subset $\{(0,0),(1,0)\}.$ This corresponds to the case where
$L_{\lambda,t}$ does not have a step in $z=1.$
Let $P_{\mbox{step}}:=P\circ\pi_{\{(0,0),(1,1)\}}^{-1}$ denote the marginal distribution
of $P$ on the subset $\{(0,0),(1,1)\}.$ This corresponds to the case where
$L_{\lambda,t}$ has a step in $z=1.$
Then, the limits in both (i) and (ii) of the above corollary are of the form
\begin{equation}
\lambda\int_{\Upsilon^2} f(y)\,P_{\mbox{step}}(dy)
+(1-\lambda)\int_{\Upsilon^2} f(y)\,P_{\mbox{flat}}(dy).
\end{equation}

\nocite{Mil88}

\section*{Acknowledgements}

This work was motivated by a special case of Corollary \ref{Unif:ex.staircase}
which originally occured during my thesis work.
I am grateful to my supervisor Hans F{\"o}llmer for the
rich education he gave me in both ergodic theory and probability theory.
It is a pleasure for me to thank Jon Aaronson
for guiding me toward a broader perspective for Lemma \ref{ErgS}.
David Steinsaltz assisted with an approximation
argument I used in the proof of Theorem \ref{UnifL}, and I thank him for this hint.
Benjamin Weiss and Yuval Peres helped me situate my results
within the field of ergodic theory;
I am grateful to them for the interest they showed in my work.
My thanks are also due to Hans Crauel, Manfred Denker,
Didier Piau, Michael Scheutzow and Masha Saprykina for helpful discussions
on different aspects of ergodic theory.  
I would like to express my gratitude to Bill MacKillop for supporting me in 
manifold ways during the time I've been working at Queen's University.



\smallskip

\noindent
{\em
Julia Brettschneider\\[1.5mm]
Department of Statistics, University of Warwick, Coventry, CV4 7AL, UK\\[1.5mm]   
Department of Community Health $\&$ Epidemiology and Cancer Research Institute Division of \\
Cancer Care $\&$ Epidemiology, Queen's University, Ontario, K7L 3N6, Canada  \\[1.5mm]    
{\tt  julia.brettschneider@warwick.ac.uk}               
}

\end{document}